\documentclass[11pt]{amsart}
\usepackage[utf8]{inputenc}
\usepackage[english]{babel}
\usepackage{amsthm,amsmath,amssymb}
\usepackage{textcomp}
\usepackage{amssymb}
\usepackage{amsthm}
\usepackage{graphicx}

\usepackage{todonotes}
\usepackage{mathtools}
\usepackage{url}
\usepackage[T1]{fontenc}
\providecommand{\Q}{\mathbb{Q}}
\newcommand{\prob}{\mathbb{P}}	
\newcommand{\expect}{\mathbb{E}}	
\newcommand{\ds}{\displaystyle}
\theoremstyle{plain}
\newtheorem{theorem}{Theorem}[section]

\newtheorem{conjecture}[theorem]{Conjecture}

\newtheorem{definition}[theorem]{Definition}

\theoremstyle{definition}
\theoremstyle{remark}

\title{Maximal persistence in random clique complexes}
\author{Ayat Ababneh}
\address{Department of Mathematics, The University of Jordan, Amman 11942, Jordan}
\author{Matthew Kahle}
\address{Department of Mathematics, The Ohio State University, Columbus, OH, USA}
\thanks{Both authors gratefully acknowledge the support of NSF-DMS \#2005630.}
\date{\today}

\begin{document}

\maketitle

\begin{abstract}
We study the persistent homology of an Erd\H{o}s--R\'enyi random clique complex filtration on $n$ vertices. Here, each edge $e$ appears at a time $p_e \in [0,1]$ chosen uniform randomly in the interval, and the \emph{persistence} of a cycle $\sigma$ is defined as $p_2 / p_1$, where $p_1$ and $p_2$ are the birth and death times of the cycle respectively. We show that for fixed $k \ge 1$, with high probability the maximal persistence of a $k$-cycle is of order roughly $n^{1/k(k+1)}$. These results are in sharp contrast with the random geometric setting where earlier work by Bobrowski, Kahle, and Skraba shows that for random \v{C}ech and Vietoris--Rips filtrations, the maximal persistence of a $k$-cycle is much smaller, of order $\left(\log n / \log \log n \right)^{1/k}$.
\end{abstract}

\section{Introduction}
Recently, the topology of random simplicial complexes has been an active area of study --- see, for example, the surveys \cite{kahle2017random} and \cite{BK18}. This study has had applications in topological data analysis, including in neuroscience \cite{GPCI}. 
We will assume that the reader is familiar with the notions of persistent homology and a persistence diagram \cite{ELZ01}. In topological inference, one sometimes considers points far from the diagonal in the persistence diagram to be representing ``signal'' and points near the diagonal as representing ``noise''. 

With this in mind, Bobrowski, Kahle, and Skraba studied maximally persistent cycles in random geometric complexes in \cite{BKS17}. They showed that with high probability the maximal persistence of a $k$-dimensional cycle in a random geometric complex in $\mathbb{R}^d$ is $\asymp \left(\log n / \log \log n \right)^{1/k}$.
We write $f \asymp g$ if $f$ and $g$ grow at the same rate in the sense that there exist constants $c_1, c_2 > 0$ such that $c_1 f(n) \le g(n) \le c_2 f(n)$ for all large enough $n$.  Both the Vietoris--Rips and \v{C}ech filtrations have an underlying parameter $r$. Persistence of a cycle is measured multiplicatively as $r_2 / r_1$ where $r_1$ and $r_2$ are the birth and death radius.

Our main result is that for fixed $k \ge 1$, the maximal persistence of $k$-cycles in a random clique complex filtration is of order roughly $n^{1/k(k+1)}$. A precise statement is given in Theorem \ref{thm:main}. Here we measure the persistence of a cycle as $p_2 / p_1$, where $p_1$ and $p_2$ are the birth and death edge probabilities, respectively. 

The comparison between the Erd\H{o}s--R\'enyi and random geometric settings may be more apparent if we renormalize so that the persistence of the associated filtrations can be measured on the same scale. One natural way to try to do this is to reconsider maximal persistence in random geometric complexes, using instead birth and death edge probability rather than radius.

The edge probability $P$ in a random geometric complex is of order $P \asymp r^d$. So if
\[
r_2 / r_1  \asymp \left( \log n / \log \log n \right)^{1/k},
\]

then
\[
P_2 / P_1  \asymp  \left( \log n / \log \log n \right)^{d/k}.
\]

\medskip

As long as $d$ is fixed, $\left( \log n / \log \log n \right)^{d/k}$ is still much smaller than the maximal persistence of cycles in the random clique complex. However, this parameterization makes it clear that if $d$ grows, we expect cycles to persist for longer. It is known that random geometric graphs in dimension growing quickly enough converge in total variation distance to Erd\H{o}s-R\'enyi random graphs, and this connection has been further explored and quantified in a number of recent papers --- see, for example \cite{BDER16, BBN20, AndrewElliot}. From this point of view, our main result can be seen as a ``curse of dimensionality'' for topological inference---as the ambient dimension gets bigger, noisy cycles persist for much longer.

\section{Topology of random clique complexes}

In this section, we review the definition of the random clique complex, and briefly survey the literature on topology of random clique complexes. \\

The following random graph is sometimes called the Erd\H{o}s--R\'{e}nyi model.
\begin{definition}\label{def:random graph}
For $n \ge 1$ and $p \in [0,1]$, $G(n,p)$ is the probability space of all simplicial graphs on $n$ vertices where every edge is included with probability $p$, jointly independently.
\end{definition}

We use the notation $G\sim G(n,p)$ to indicate that a graph $G$ is chosen according to this distribution. We say that $G$ has a given property \emph{with high probability (w.h.p.)} if the probability that $G$ has the property tends to one as $n \to \infty$.

We define $X(n,p)$ to be the clique complex of the Erd\H{o}s--R\'enyi random graph $G(n,p)$. We write $X\sim X(n,p)$ to indicate that $X$ is a random simplicial complex chosen according to this distribution.

In ~\cite{Kahle2009}, Kahle studied the topology of random clique complexes, and the following theorem roughly identified the threshold for homology to appear.

% \begin{theorem}\label{thm:Hom.Vanishing}
% If $p=n^\alpha$ with $\alpha<-\frac{1}{k}$ then $\widetilde{H}_k\left(X(n,p), \mathbb{Z}\right)=0.$
% \end{theorem}
\begin{theorem}\label{thm:Threshold} Let $X \sim X(n,p)$ be the random clique complex, and assume $k \ge 1$ is fixed.
\begin{enumerate}
    \item If $p \le n^{-\alpha}$ with $\alpha > \frac{1}{k}$, then w.h.p.\ $H_k(X)=0$. On the other hand,
    \item If $n^{-\frac{1}{k}} \ll p \ll n^{-\frac{1}{k+1}}$ then w.h.p.\ $H_k(X)\ne 0$.
\end{enumerate}
\end{theorem}

We use the notation $f \ll g$ to indicate $\lim_{n \to \infty} f/g = 0$. Theorem \ref{thm:Threshold} shows that the threshold for the appearance of $k$th homology is roughly $n^{-1/k}$. The following is the main result of a later paper, \cite{Kahle14}, showing that the threshold for the vanishing of $k$th homology is roughly $n^{-1/(k+1)}$.

\begin{theorem}\label{thm:cohomology}
Let $k\geq1$ and $\epsilon>0$ be fixed, and $X\sim X(n,p)$. If 
\[
p\geq\left(\frac{(\frac{k}{2}+1+\epsilon)\log n}{n}\right)^{1/(k+1)}
\]
then w.h.p. $H^k(X,\Q)=0$.
\end{theorem}

Theorem \ref{thm:cohomology} describes a sharp threshold for cohomology to vanish, in the same spirit as in Linial and Meshulam's work \cite{LM06}. These theorems all describe high-dimensional cohomological generalizations of the Erd\H{o}s--R\'enyi theorem on the threshold for connectivity of $G(n,p)$ \cite{ErRe59}. By universal coefficients, since we are working over a field, the results hold for homology with $\Q$ coefficients as well as cohomology.

Malen gave a topological strengthening of part (1) of Theorem \ref{thm:Threshold} in \cite{Malen2019}.

\begin{theorem}[Malen, 2019] Let $k \ge 1$ be fixed and $X \sim X(n,p)$. If $p \le n^{-\alpha}$ with $\alpha > \frac{1}{k}$, then w.h.p.\ $X$ collapses onto a subcomplex of dimension at most $k-1$.
\end{theorem}

This implies, in particular, that $H_{k-1}(X)$ is torsion-free, so this represents an important step toward the ``bouquet-of-spheres conjecture'' described in \cite{Kahle2009} and \cite{Kahle14}.

Newman recently refined Malen's collapsing argument to give a probabilistic refinement \cite{Newman2021}.

\begin{theorem} [Newman, 2021] \label{thm:birth}
Let $k \ge 1$ be fixed and $X \sim X(n,p)$.
If 
\[
p \ll n^{-1/k}
\]
then w.h.p.\ $X$ collapses onto a subcomplex of dimension at most $k-1$.
\end{theorem}

In summary, earlier results show that there is one threshold where homology is born for the first time, when $p \approx n^{-1/k}$, and another where homology dies for the last time, when $p \approx n^{-1/(k+1)}$. Our main result is that there exist cycles that persist for nearly the entire interval of nontrivial homology.

\section{The second moment method}

We briefly review the second moment method, i.e.\ the use of Chebyshev's inequality, which is our main probabilistic tool.

\begin{theorem}[Chebyshev's Inequality]\label{thm:chebyshev}
For any $\lambda>0$, 
\[
\prob\left(\lvert X-\mu\rvert\geq\lambda\sigma\right)\leq\frac{1}{\lambda^2}.
\]
Where $\mu$ is the expectation and $\sigma^2$ is the variance.
\end{theorem}
The variance is defined by
\[
\sigma^2 = \text{Var}(X)=\expect\left(X-\expect(X)^2\right).
\]
If $X$ can be written as a sum of indicator random variables $X=\ds \sum_i X_i$, then the following is easy to derive and its proof appears, for example, in Chapter 4 of Alon and Spencer's book \cite{AS16}.
\[
\begin{array}{lll}
\text{Var}(X)&=&\ds \sum_{i} \text{Var}(X_i)+\sum_{i \neq j} \text{Cov}(X_i,X_j)
\\\\
&\leq& \expect(X)+\ds \sum_{i \neq j} \text{Cov}(X_i,X_j).
\end{array}
\]

It follows from Chebyshev's inequality \ref{thm:chebyshev} that if $\expect(X)\rightarrow\infty$ and
\[
\ds \sum_{i \neq j} \mbox{Cov}(X_i,X_j) = o\left(\expect(X)^2\right),
\]
then $X > 0$ w.h.p. In fact, $X \sim \expect(X)$ w.h.p., meaning that $X / \expect(X) \to 1$ in probability.

Finally, we note that if are  $X_i, X_j$ are indicator random variables for events $A_i, A_j$, we have that
\[
\text{Cov}(X_i,X_j) = \expect(X_iX_j)-\expect(X_i)(X_j)= \prob(A_i \wedge A_j) - \prob(A_i)\prob(A_j).
\]
Here $A_i \wedge A_j$ denotes the event that both $A_i$ and $A_j$ occur. 

\section{Main result and proof}

We consider the random graph $G(n,p)$ as a stochastic process, as follows. Consider the random filtration of the complete graph $K_n$ where each edge $e$ appears at time $p_e$, chosen uniform randomly in the interval $[0,1]$. Similarly, the random clique complex $X(n,p)$ is a random filtration of the simplex on $n$ vertices $\Delta_n$.

We assume the reader is familiar with persistent diagrams \cite{ELZ01,CsEG07}. A point $(x,y)$ in the persistence diagram for $H_k$ with $\Q$ coefficients and $k \ge 1$ represents a $k$-dimensional cycle with birth time $x$ and death time $y$. We measure the persistence of that cycle multiplicatively, as $y/x$. Define
\[
M_k (n) = \max \{ y/x \},
\]
where the maximum is taken over all points in the persistence diagram for homology in degree $k$.

An equivalent definition is the following. Consider the natural inclusion map $i: X(n,p_1) \hookrightarrow X(n,p_2)$, where $0 \le p_1 \le p_2 \le 1$. For every $k \ge 1$, there is an induced map on homology $i_*: H_k(X(n,p_1)) \to H_k(X(n,p_2))$. Then
\[
M_k(n)= \max \left\{ p_2 / p_1 \mid i_*: H_k(X(n,p_1)) \to H_k(X(n,p_2)) \mbox{ is nontrivial} \right\},
\]
where the maximum is taken as $p_1$ and $p_2$ range over all values with $0 \le p_1 \le p_2 \le 1$.

\bigskip

Our main result is the following. \\

\begin{theorem} \label{thm:main}
For fixed $k\geq1$ and $\epsilon > 0$,
\[
n^{-1/k(k+1)-\epsilon} \le M_k(n) \le  n^{-1/k(k+1)+\epsilon},
\]
with high probability.

Equivalently, if
\[
\widetilde{M}_k(n) = \frac{ \log {M_k(n)}}{\log{n}},
\]
then $\widetilde{M}_k(n)$ converges in probability to 
$1/k(k+1).$ 
\end{theorem}

Our results are actually slightly sharper then this. A slightly sharper statement is the following. 

Suppose that 
\[
L_k(n) \ll n^{1/k(k+1)}
\]and
\[
U_k(n) \gg n^{1/k(k+1)} (\log n)^{1/(k+1)}.
\]
Then w.h.p.
\[ L_k(n) \le M_k(n) \le U_k(n). \]
So then our results are sharp, up to a small power of $\log n$.\\

Since Theorem \ref{thm:cohomology} is only known to hold for $\Q$ coefficients (or any field of characteristic zero), our upper bounds on maximal persistence only apply for persistent homology with $\Q$ coefficients. For lower bounds, we use the second moment method to prove the existence of cycles that persist nearly as long as possible, and these lower bounds hold for arbitrary field coefficients. \\

\begin{proof}
First we prove an upper bound on $M_k(n)$.

Suppose that $p_1 = o \left(n^{-1/k} \right)$. By Theorem \ref{thm:birth}, we have that $H_k(X) = 0$ w.h.p. Now let $\epsilon > 0$, and suppose that 
\[
p_2 \geq\left(\frac{(\frac{k}{2}+1+\epsilon)\log n}{n}\right)^{1/(k+1)}.
\]
By Theorem \ref{thm:cohomology}, we have that $H_k(X,Q) = 0$. So we have that
\[
M_k(n) \le p_2 / p_1.
\]

Set
\[
f_k(n) = n^{1/k(k+1)} \left( \log{n} \right)^{1/(k+1)},
\]
and let $U_k(n)$ be any function such that $U_k(n) \gg f_k(n)$. We have showed that
\[
M_k(n) \le U_k(n).
\]

\bigskip

Most of our work is in proving a lower bound for $M_k(n)$. We focus our attention on a particular type of nontrivial $k$-cycle, namely simplicial spheres which are combinatorially isomorphic to cross-polytope boundaries.

In the following, let $Y$ and $Z$ denote distinct subsets of $2k+2$ vertices. That is, we write $[n]:= \{ 1, 2, \dots, n\}$ and suppose that $Y, Z \subseteq [n]$ with $|Y|=|Z|=2k+2$. A notation we can use for this is 
\[
Y, Z \in \binom{[n]}{2k+2}.
\]

Suppose that $Y= \{ u_1, \dots u_{k+1} \} \cup \{v_1, \dots ,v_{k+1}\}$, where $u_1 < \dots u_{k+1} < v_1 \dots < v_{k+1}$. Recall that every vertex is an element of $[n]$, so they come with a natural ordering. We use $x \sim y$ and $x \not\sim y$ to denote adjacency and non-adjacency of vertices $x$ and $y$. For any choice of $0 \le p_1 \le p_2 \le 1$, we say that $Y$ is a \emph{$(p_1,p_2)$ special persistent cycle} in the random clique complex filtration if 
\begin{enumerate}
    \item $u_i \sim u_j$, $v_i \sim v_j$, and $u_i \sim v_j$ for every $i \neq j$ at time $p_1$,
    \item $u _ i \not\sim v_i$ for every $i$ at time $p_2$, and
    \item $\{ u_1, \dots u_{k+1}\}$ have no common neighbors outside of vertex set $Y$ at time $p_2$.
\end{enumerate}
Condition (1) implies that $Y$ spans a $k$-dimensional cycle at time $p_1$, namely a cycle that is combinatorially equivalent to the boundary of $k+1$-dimensional cross-polytope. Conditions (2) and (3) imply that $Y$ is still not a boundary of anything, even at time $p_2$. So then it is not only a nontrivial cycle at time $p_1$, but it persists at least until time $p_2$. Note that condition (2) already implies that $\{ u_1, \dots u_{k+1}\}$ have no common neighbor within vertex set $Y$. So condition (3) implies that they have no common neighbor at all, and then $\{ u_1, \dots u_{k+1} \}$ is a maximal $k$-dimensional face. 

Let $N_k=N_k(p_1,p_2)$ be the number of $(p_1,p_2)$ special persistent cycles. We want to show that $\prob\left(N_k>0\right)\rightarrow 1$, which in turn will imply that $M_k(n) > p_2 / p_1$ with high probability. In the following, we will assume whenever necessary that
\[
n^{-1/k} \ll p_1 \le p_2 \ll n^{-1/(k+1)}.
\]
In particular, we assume that $n p_1^k \to \infty$ and $n p_2^{k+1} \to 0$.

Let $A_Y$ be the event that the set of vertices in $Y$ form a $(p_1,p_2)$ special persistent cycle, and let $I_Y$ be its indicator random variable for this event. Then we can write 
\[
N_k=\ds\sum_{Y\in \binom{[n]}{2k+2}} I_Y,
\]
where the sum is taken over all subsets $Y \subseteq [n]$ of size $|Y|=2k+2$.

By edge independence, the probability of condition (1) is $p_1^{2k(k+1)}$, and the probability of condition (2) is $\left(1-p_2\right)^{k+1}$, the probability of condition (3) is $\left(1-p_2^{k+1}\right)^{n-2k-2}$. Moreover, these events are independent since they involve disjoint sets of edges. So we have
\[
\expect(I_Y)=\prob(A_Y)=p_1^{2k(k+1)}\left(1-p_2\right)^{k+1}\left(1-p_2^{k+1}\right)^{n-2k-2}.
\]
By linearity of expectation,
\[
\begin{split}
\expect(N_k)&=\sum \expect(I_Y) \\
&=\ds \binom{n}{2k+2}p_1^{2k(k+1)}\left(1-p_2\right)^{k+1}\left(1-p_2^{k+1}\right)^{n-2k-2}\\
&= \ds\frac{n^{2k+2}}{\left(2k+2\right)!}p_1^{2k(k+1)} \left( 1-o(1) \right),
\end{split}
\]
since $np_2^{k+1} \to 0$.
Since we also assume that $n p_1^{k}\to \infty$, we have $\expect(N_k)\rightarrow\infty$. By Chebyshev's inequality, if we show that $\text{Var}(N_k)=o(\expect(N_k)^2)$, then $N_k > 0$ w.h.p.

We have the standard inequality
\[
\text{Var}(N_k)\leq \expect(N_k)+\sum_{Y \neq Z}\text{Cov}(I_Y,I_Z).
\]

We recall that 
\[
\text{Cov}(I_Y,I_Z) = \prob(A_Y \wedge A_Z) - \prob(A_Y) \prob(A_Z).
\]

\bigskip

We always have
\[
\prob(A_Y)\prob(A_Z)=p_1^{4k(k+1)}\left(1-p_2\right)^{2k+2}\left(1-p_2^{k+1}\right)^{2(n-2k-2)}, 
\]
and we note the simpler estimate 
\[
\prob(A_Y)\prob(A_Z)=p_1^{4k(k+1)} \left( 1- o(1) \right)
\]
since $k \ge 1$ is fixed and $np_2^{k+1} \to 0$.

Let $m:=\rvert Y\cap Z\lvert$. In estimating $\prob(A_Y \wedge A_Z)$, we consider cases depending on the value of $m$. \\

\bigskip

{\bf Case I:}

First, consider $m=0$. It might be tempting to believe that in the case that $Y \cap Z = \emptyset$, $A_Y$ and $A_Z$ are independent sets, so the covariance is zero, but unfortunately this is not the case. Conditions (1) and (2) for a $(p_1,p_2)$ special persistent cycle only depend on adjacency between vertices within the $(2k+2)$-set, but condition (3) depends on connections with the rest of the graph and these are not independent.

Nevertheless, we still have in this case
\[
\prob\left(A_Y \wedge A_Z\right)= p_1^{4k(k+1)}\left(1-o(1) \right),
\]
as follows.

The term $p_1^{4k(k+1)}$ is the probability of condition (1) holding for both vertex sets $Y$ and $Z$. So this is also an upper bound on the probability of conditions (1), (2), and (3) holding for both vertex sets. For a lower bound on $\prob\left(A_Y \wedge A_Z\right)$, we consider a slightly smaller event, slightly simpler but whose probability is of the same order of magnitude.

Let $Y= \{ u_1, \dots u_{k+1} \} \cup \{v_1, \dots ,v_{k+1}\}$, where $u_1 < \dots u_{k+1} < v_1 \dots < v_{k+1}$, as before. Similarly, let $Z= \{ u'_1, \dots u'_{k+1} \} \cup \{v'_1, \dots ,v'_{k+1}\}$, where $u'_1 < \dots u'_{k+1} < v'_1 \dots < v'_{k+1}$.

The event $A^*_{YZ}$ is defined as follows. 

\begin{enumerate}
    \item We have $u_i \sim u_j$, $v_i \sim v_j$, and $u_i \sim v_j$, $u'_i \sim u'_j$, $v'_i \sim v'_j$, and $u'_i \sim v'_j$ for every $i \neq j$ at time $p_1$. That is, condition (1) holds for both $Y$ and $Z$. Some edges may be listed more than once if $Y$ and $Z$ overlap. This does not happen when $m=0$ but these are the cases we consider below.
    \item Besides the edges that appear in the previous condition other edges occur between vertices in vertex set $Y \cup Z$, at time $p_2$. This happens with probability $1-\mathcal{O}(p_2) = 1-o(1)$.
    \item Neither $\{ u_1, \dots u_{k+1} \}$ nor $\{ u'_1, \dots u'_{k+1} \} $ has any mutual neighbors outside of vertex set $Y \cup Z$, at time $p_2$. The probability of this condition being satisfied can be bounded below by a union bound by $1-2np_2^{k+1}$, which is again $1-o(1)$ since $np_2^{k+1} \to 0$.
\end{enumerate}

It is clear that $A^*_{YZ}$ imples $A_Y \wedge A_Z$. Indeed, condition (1) is the same, condition (2) of $A^*_{YZ}$ implies condition (2) of $A_Y \wedge A_Z$, and conditions (2) and (3) of $A^*_{YZ}$ together imply condition (3) of $A_Y \wedge A_Z$. So putting it all together, we have that
$\prob( A^*_{YZ}) \ge p_1^{4k(k+1)} \left( 1-o(1) \right)$.

Then
\[
p_1^{4k(k+1)} \ge \prob (A_Y \wedge A_Z ) \ge \prob ( A^*_{YZ} ) \ge p_1^{4k(k+1)} \left (1 - o(1) \right),
\]
and 
\[ 
\prob (A_Y \wedge A_Z ) = p_1^{4k(k+1)}\left (1 - o(1) \right),
\]
as desired. \\

So then
\[
\prob \left(A_Y \wedge A_Z\right)-\prob (A_Y)\prob (A_Z)= o\left( p_1^{4k(k+1)}\right).
\]
Since the number of pairs $Y,Z$ is bounded by $n^{4k+4}$ we have that the total contribution to the variance, $S_0$, is bounded by
\[
S_0=o\left(n^{4k+4}p_1^{4k(k+1)}\right).
\]
Comparing this to 
\[
\expect(N_k)^2=\ds\binom{n}{2k+2}^2 p_1^{4k(k+1)}\left(1-o(1)\right)
\]
we see that
\[
S_0=o \left( \expect(N_k)^2 \right).
\]

\bigskip

{\bf Case II:}

An essentially identical calculation shows that when $m=1$, we have
\[
\prob \left(A_Y \wedge A_Z\right)= p_1^{4k(k+1)} \left(1-o(1) \right).
\]
So in this case we have again
\[
\prob \left(A_Y \wedge A_Z\right)-\prob (A_Y)\prob (A_Z)= o\left( p_1^{4k(k+1)}\right).
\]
Hence, the total contribution to the variance, $S_1$, is
\[
S_1=o \left(n^{4k+3}p_1^{4k(k+1)}\right)
\]
and then,
\[
S_1 / \expect(N_k)^2 = o \left( n^{-1} \right)
\]
and in particular $S_1=o\left(\expect(N_k)^2\right)$.\\

{\bf Case III:}\\

When $2\leq m\leq 2k+1$, we consider two sub-cases. The first subcase is that events $A_Y$ and $A_Z$ are not compatible in the sense that they cannot both occur due to the ways in which $Y$ and $Z$ overlap. This happens if for a certain pair of vertices $u,v \in Y \cap Z$, $u,v$ are required to be adjacent in one of $Y, Z$ and non-adjacent in the other. In this subcase, we have
\[
\prob ( A_Y \wedge A_Z) = 0,
\]
so
\[
\prob ( A_Y \wedge A_Z) - \prob (A_Y) \prob (A_Z) = -p^{4k(k+1)} \left( 1-o(1) \right) \le 0.
\]

The second subcase is that the events $A_Y$ and $A_Z$ are compatible, in the sense that they could possibly both happen. In this case, let $j$ denote the number of pairs in $Y\cap Z$ that are forced to be non-adjacent in $A_Y \cap A_Z$. Then the same argument as in Case I shows that 
\[
\prob \left( A_Y \wedge A_Z\right) \ge \prob \left( A^*_{YZ} \right) = p_1^{4k(k+1)-{\binom{m}{2}+j}}\left(1-o(1) \right).
\]

So
\[
\begin{array}{lll}
\ds\frac{\prob (A_Y)\prob (A_Z)}{\prob (A_Y \wedge A_Z)}&\leq& p_1^{{\binom{m}{2}}-j}\left(1-o(1)\right)
\end{array}
\]
Since $np_1^{k}\rightarrow\infty,$ as $n\rightarrow\infty$, we get 
\[
\ds\frac{\prob (A_Y)\prob (A_Z)}{\prob (A_Y \wedge A_Z)}\rightarrow 0.
\]
So,
\[
\prob \left(A_Y \wedge A_Z\right)-\prob (A_Y)\prob (A_Z)=\left(1-o(1)\right)\prob \left(A_Y \wedge A_Z\right).
\]
The total contribution $S_m$ of a pair of events $A_Y$ and $A_Z$ with $Y\cap Z =m$ to the variance is then bounded by
\[
S_m \le n^{4k+4-m}p_1^{4k(k+1)-{\binom{m}{2}}+j}\left( 1+ o(1) \right).
\]
Comparing this to
\[
\expect(N_k)^2=\ds \binom{n}{2k+2}^2 p_1^{4k(k+1)}\left(1-o(1)\right).
\]
we get
\[
S_m/\expect(N_k)^2=O\left(n^{-m}p_1^{-{\binom{m}{2}}+j}\right).
\]
We have
\[
n^{-m}p_1^{-{\binom{m}{2}}+j}=\left(np_1^{\frac{m-1}{2}}\right)^{-m} p_1^j.
\]
We are assuming  that $n p_1^{k} \rightarrow \infty$. Since $m \le 2k+1$, we have $k \ge (m-1)/2$. Then
\[
\left(np_1^{\frac{m-1}{2}}\right)^{-m}\rightarrow 0,
\]
$p_1^j \to 0$, and $S_m=o\left(\expect(N_k)^2\right)$. \\

Summing the inequalities from the different cases, we conclude that
\[
\sum_{Y \neq Z}\text{Cov}(I_Y,I_Z) = \sum_{m=0}^{2k+1} S_m = o\left(\expect(N_k)^2\right), 
\]
since $S_m =o\left(\expect(N_k)^2\right)$ for each $m$ and $k$ is fixed.

We conclude that as long as
\[
n^{-1/k} \ll p_1 \le p_2 \ll n^{-1/(k+1)},
\]
then $N_k > 0 $ with high probability.

It follows that if $L_k(n) \ll n^{1/k(k+1)}$ then w.h.p. $M_k(n) \ge L_k(n)$, as desired.

\end{proof}
\section{Future directions}

Recall that we earlier defined
\[
f_k(n) = n^{1/k(k+1)} \left( \log{n} \right)^{1/(k+1)}.
\]
We believe that the $M_k(n)$ is likely of order roughly $f_k(n)$, in the following sense.

Let $\omega(n)$ be any function that tends to infinity with $n$.
We showed in the proof of Theorem \ref{thm:main} that 
\[
M_k(n) \le f_k(n) \omega(n).
\]
We believe that an analogous lower bound should hold.

\begin{conjecture}
Let $M_k(n)$ denote the maximal persistence over all $k$-dimensional cycles in $X(n,p)$.  Then
\[
\frac{f_k(n)}{\omega(n)}  \le M_k(n)
\]
with high probability.
\end{conjecture}

The following kind of limit theorem would provide precise answers to questions like, ``Given a prior of this kind of distribution, what is the probability $P(\lambda)$ that there exists a cycle of persistence greater than $\lambda$?''
\begin{conjecture} \label{conj:dist}
Let $M_k(n)$ denote the maximal persistence over all $k$-dimensional cycles in $X(n,p)$. Then
\[
\frac{M_k(n)}{f_k(n)}
\]
converges in law to a limiting distribution supported on an interval $[\lambda_k, \infty)$ for some $\lambda_k > 0$.
\end{conjecture}

\begin{figure}
    \centering
    \includegraphics[width=2.25in]{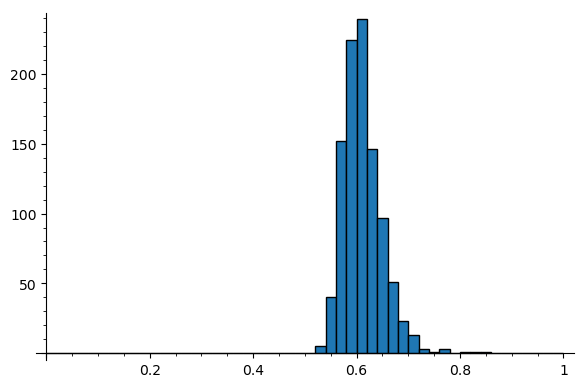}
    \includegraphics[width=2.25in]{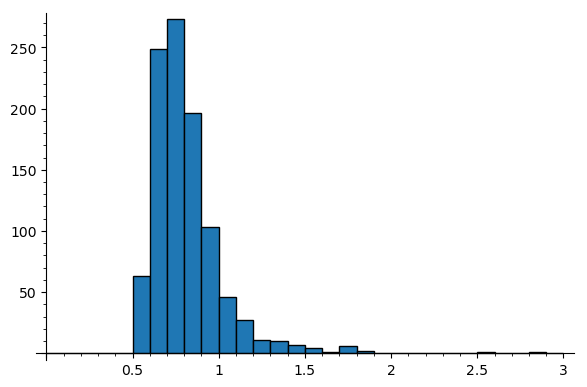}
    \caption{Maximally persistent $1$-cycles. On the left, a histogram for $\log M_1(n) / \log n$. We prove that this converges to $1/2$ as $n \to \infty$. On the right, a histogram for $M_1(n) / f_1(n)$. Both these figures are based on 1000 samples on $n=250$ vertices.}
    \label{fig:experiments1}
\end{figure}

\begin{figure}
    \centering
    \includegraphics[width=2.25in]{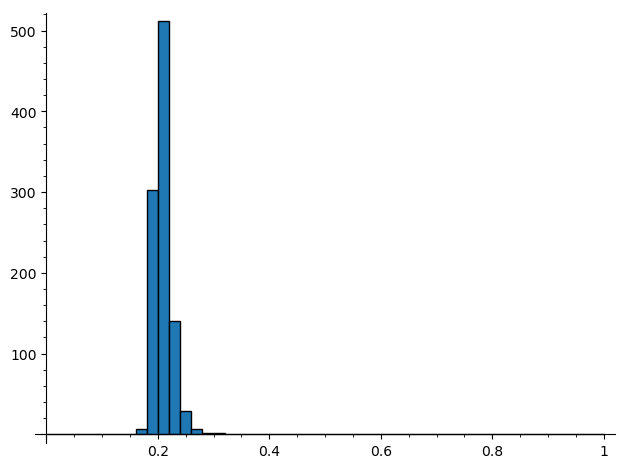}
    \includegraphics[width=2.25in]{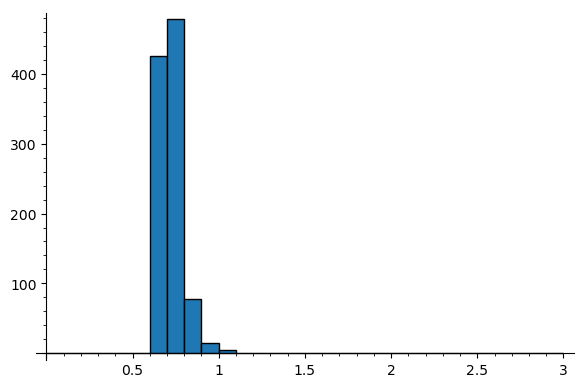}
    \caption{Maximally persistent $2$-cycles. On the left, a histogram for $\log M_2(n) /  \log n$. We prove that this converges to $1/6$ as $n \to \infty$. On the right, a histogram for $M_2(n) / f_2(n)$. Both these figures are based on 1000 samples on $n=150$ vertices.}
    \label{fig:experiments2}
\end{figure}

See Figures \ref{fig:experiments1} and \ref{fig:experiments2} for some numerical experiments illustrating these conjectures. These experiments were computed with the aid of Ulrich Bauer's software Ripser \cite{Bauer2021Ripser}.

It also seems natural to study more about the ``rank invariant'' of a random clique complex filtration. That is, given $k \ge 1$, $p_1$, and $p_2$, how large do we expect the rank of the map 
$i_* : H_k \left( X(n,p_1) \right) \to H_k \left( X(n,p_2) \right)$ to be?

\begin{conjecture}
Suppose that $k \ge 1$ is fixed, and
\[
n^{-1/k} \ll p_1 \le p_2 \ll n^{-1/(k+1)}.
\]
If $i: X(n,p_1) \to X(n,p_2)$ is the inclusion map, and
\[
i_* : H_k \left( X\left(n,p_1\right) \right) \to H_k\left( X \left(n,p_2 \right) \right)
\]
is the induced map on homology, then
\[
\mbox{rank}(i_*) = \left( 1-o(1) \right)  \binom{n}{k+1} p_1^{\binom{k+1}{2}}.
\]
\end{conjecture}
In \cite{Kahle2009}, it is shown that 
\[ \dim H_k(X(n,p_1),\Q) =  \left( 1-o(1) \right)  \binom{n}{k+1} p_1^{\binom{k+1}{2}},
\]
so this conjecture is that almost all of the homology persists for as long as possible. 

Bobrowski and Skraba study limiting distributions for maximal persistence in their recent preprint \cite{BS22}. They describe experimental evidence that there is a universal distribution for maximal persistence over a wide class of models, including random \v{C}ech and Vietoris--Rips complexes. It is not clear whether we should expect the random clique complex filtration studied here to be in the same conjectural universality class.

\bibliographystyle{plain}
\bibliography{references}

\end{document}